\documentclass{article}
\usepackage{amsmath,amsthm,amsfonts}
\usepackage{amssymb}
\usepackage{mathrsfs}
\usepackage{euscript}
\newtheorem{theorem}{Theorem}[section]
\newtheorem{definition}{Definition}[section]

\numberwithin{equation}{section}
\title{Characterization of the generalized Chebyshev-type polynomials of first kind}
\author{Mohammad A. AlQudah\\
  \small Department of Mathematics, Northwood University\\ 
	\small 4000 Whiting Dr, Midland, MI 48640 USA\\
    \small email: alqudahm@northwood.edu}

\begin{document}
\maketitle
\begin{abstract} Orthogonal polynomials have very useful properties in the solution of mathematical problems, so recent years have seen a great deal in the
 field of approximation theory using orthogonal polynomials. In this paper, we characterize the generalized Chebyshev-type polynomials of the 
first kind  $\mathscr{T}_{n}^{(M,N)}(x),$  then we provide a closed form of the constructed polynomials in term of the Bernstein polynomials $B_{k}^{n}(x).$
We conclude the paper with some results on the integration of the weighted generalized Chebyshev-type with the
 Bernstein polynomials.
\end{abstract} 



\section{Introduction and Background}
Approximation is essential to many numerical techniques, since it is possible to approximate arbitrary continuous function by a polynomial, at the same time polynomials can be represented in many different bases such as monomial and Bernstein basis.
\subsection{Univariate Chebyshev polynomials} Chebyshev polynomial of first kind (Chebyshev-I) of degree $n\geq 0$ in $x$ is defined as 
$T_{n}(x)=\cos(n \arccos x),$ $x\in[-1,1].$ The polynomials $T_{n}(x)$ of degree $n$ are orthogonal polynomials, except for a constant factor, with respect to the weight function
$\mathrm{W}(x)=\frac{1}{\sqrt{1-x^{2}}}.$ They are a special case of Jacobi polynomials $P_{n}^{(\alpha,\beta)}(x)$ and the interrelation is given by
\begin{equation}\label{cheb-Jac-Rel}T_{n}(x)=\binom{n-\frac{1}{2}}{n}^{-1} P_{n}^{(-\frac{1}{2},-\frac{1}{2})}(x).\end{equation}

Authors are not uniform on notations, and for the convenience we recall few representations of univariate Chebyshev-I polynomials. It is worth pointing out that 
univariate classical orthogonal polynomials are traditional defined on $[-1,1],$ however, it is more convenient to use $[0,1].$
The univariate Chebyshev-I polynomials of degree $n$ in $x$ can be written as \cite{KoekoekB,Szego}. 
\begin{equation}\label{Chebyshev-I-univ}
T_{n}(x):=\frac{(2n)!!}{(2n-1)!!}\sum_{k=0}^{n}\binom{n-\frac{1}{2}}{n-k}\binom{n-\frac{1}{2}}{k}\left(\frac{x-1}{2}\right)^{k}\left(\frac{x+1}{2}\right)^{n-k},
\end{equation}
which it can be transformed in terms of Bernstein basis on $x\in[0,1]$ as
\begin{equation}\label{Chebyshev-I-basis fromat}
T_{n}(2x-1):=\frac{2^{2n}(n!)^{2}}{(2n)!}\sum_{k=0}^{n}(-1)^{n+1}\frac{\binom{n-\frac{1}{2}}{k}\binom{n-\frac{1}{2}}{n-k}}{\binom{n}{k}}B_{k}^{n}(x),
\end{equation}
where $B_{k}^{n}(x)$ are the Bernstein polynomials of degree $n,$ $x\in [0,1], k=0,1,\dots, n,$ are defined by
\begin{equation}\label{binomial} B_{k}^{n}(x)=\binom {n}{k} x^{k}(1-x)^{n-k}, \end{equation}
where the binomial coefficients $\binom{n}{k}=\frac{n!}{k!(n-k)!},\hspace{.1in}k=0,1,\dots,n.$

Note that the double factorial of an integer $m$ is a generalization of the usual factorial $m!,$ the notation $m!!$ appears not to be widely known and defined as
\begin{equation}\label{douvle-factrrial}
\begin{aligned}
m!!&=(m)(m-2)(m-4)\dots(4)(2) \hspace{.72in} \text{if $m$ is even}\\
(2m-1)!!&=(2m-1)(2m-3)(2m-5)\dots(3)(1) \hspace{.2in}  \text{if $m$ is odd},
\end{aligned}
\end{equation}
where $0!!=(-1)!!=1.$ 

From the definition \eqref{douvle-factrrial}, we can derive the factorial of an integer minus half as 
\begin{equation}\label{n-half}
\left(r-\frac{1}{2}\right)!=\frac{r!(2r-1)!!\sqrt{\pi}}{(2r)!!}.
\end{equation}
In addition, the Chebyshev-I polynomials satisfy the orthogonality relation \cite{Olver}
\begin{equation}\label{orth-rel}
\int_{0}^{1}(x-x^{2})^{-\frac{1}{2}}T_{n}(x)T_{m}(x)dx=\left\{\begin{array}{ll} 
0& \mbox{if } m\neq n \\
\frac{\pi}{2} & \mbox{if } m=n=0\\
\pi & \mbox{if } m=n=1,2,\dots
\end{array}\right..
\end{equation}

\section{Generalized Chebyshev-type polynomials}
In this section we characterize the generalized Chebyshev-type polynomials of first kind, $\mathscr{T}_{r}^{(M,N)}(x),$ then we provide explicit closed form as a linear combination of Bernstein polynomials $B_{i}^{r}(x).$
We conclude this section with the closed form of the integration of the weighted generalized Chebyshev-type with respect to the Bernstein polynomials.

Using the fact \eqref{cheb-Jac-Rel} and similar construction of the results in \cite{Bav-Koek2,Koornwinder}. For $M, N\geq 0,$  
we define the generalized Chebyshev-type polynomials
$\left\{\mathscr{T}_{n}^{(M,N)}(x)\right\}_{n=0}^{\infty}$ as
\begin{equation}\label{gen-Chebyshev-I}
\mathscr{T}_{n}^{(M,N)}(x)
=\frac{(2n)!}{2^{2n}(n!)^{2}}T_{n}(x)+MQ_{n}(x)+NR_{n}(x)+MNS_{n}(x),\hspace{.05in}n=0,1,2,\dots
\end{equation}
where for $n=1,2,3,\dots$
\begin{equation}\label{q-fn}
Q_{n}(x)=\frac{(2n)!}{2^{2n-1}(n!)^{3}}\left[n^{2}T_{n}(x)
-\frac{1}{2}(x-1)DT_{n}(x)\right], \hspace{.05in}n=1,2,3,\dots
\end{equation}
\begin{equation}\label{r-fn}
R_{n}(x)=\frac{(2n)!}{2^{2n-1}(n!)^{3}}\left[n^{2}T_{n}(x)-\frac{1}{2}(x+1)DT_{n}(x)\right], \hspace{.05in}n=1,2,3,\dots
\end{equation}
and
\begin{equation}\label{s-fn}
S_{n}(x)=\frac{(2n)!}{2^{2n-2}(n!)^{3}(n-1)!}[n^{2}T_{n}(x)
-x DT_{n}(x)], \hspace{.05in}n=1,2,3,\dots.
\end{equation}
By using
$(x^{2}-1)D^{2}T_{n}(x)=n^{2}T_{n}(x)-xDT_{n}(x),$
we find that
\begin{equation}\label{s-fn-useJac}
S_{n}(x)=\frac{4(2n-1)!!}{n!(n-1)!(2n)!!}(x^{2}-1)D^{2} T_{n}(x), \hspace{.1in}n=1,2,3,\dots.
\end{equation}
It is clear that $Q_{0}(x)=R_{0}(x)=S_{0}(x)=0.$ Furthermore, the generalized Chebyshev-type polynomials satisfy the symmetry relation \cite{Koornwinder},
\begin{equation*}
\mathscr{T}_{n}^{(M,N)}(x)=(-1)^{n}\mathscr{T}_{n}^{(N,M)}(-x),\hspace{.1in} n=0,1,2,\dots
\end{equation*}
which implies that for $n=0,1,2,\dots,$ $Q_{n}(x)=(-1)^{n}R_{n}(-x)$ and $S_{n}(x)=(-1)^{n}S_{n}(-x).$ From \eqref{q-fn} and \eqref{r-fn} it follows that for $n=1,2,3,\dots$
\begin{equation*}
Q_{n}(1)=\frac{2(2n-1)!!}{(n-1)!(2n)!!}T_{n}(1)\hspace{.1in} \text{and} \hspace{.1in} 
R_{n}(-1)=\frac{2(2n-1)!!}{(n-1)!(2n)!!}T_{n}(-1).
\end{equation*}
Note that \eqref{q-fn},\eqref{r-fn} and \eqref{s-fn} imply that for $n=1,2,3,\dots,$ we have 
\begin{equation}\label{q-fn-sum}
Q_{n}(x)=\sum_{k=0}^{n}\frac{(2k)!q_{k}}{2^{2k}(k!)^{2}} T_{k}(x)\hspace{.15in} \text{with} \hspace{.15in} 
q_{n}=\frac{4}{(2n-3)(n-1)!},
\end{equation}
\begin{equation}\label{r-fn-sum}
R_{n}(x)=\sum_{k=0}^{n}\frac{(2k)!r_{k}}{2^{2k}(k!)^{2}} T_{k}(x) \hspace{.15in} \text{with} \hspace{.15in} 
r_{n}=\frac{4}{(2n-3)(n-1)!},
\end{equation}
and
\begin{equation}\label{s-fn-sum}
S_{n}(x)=\sum_{k=0}^{n}\frac{(2k)!s_{k}}{2^{2k}(k!)^{2}} T_{k}(x) \hspace{.15in} \text{with} \hspace{.15in} s_{n}=\frac{4}{(n-1)!(n-2)!}.
\end{equation}

Therefore, for $M, N\geq 0$ the generalized Chebyshev-type polynomials\\ $\left\{\mathscr{T}_{n}^{(M,N)}(x)\right\}_{n=0}^{\infty}$ are orthogonal on the interval $[-1,1]$ with respect to the weight function
\begin{equation}
\frac{1}{\pi}(1-x)^{-\frac{1}{2}}(1+x)^{-\frac{1}{2}}+M\delta(x+1)+N\delta(x-1),
\end{equation}
and can be written as
\begin{equation}\label{gen-Chebyshev-I-ccomb}
\mathscr{T}_{n}^{(M,N)}(x)=\frac{(2n-1)!!}{(2n)!!}T_{n}(x)+\sum_{k=0}^{n} \frac{(2k)!\lambda_{k}}{2^{2k}(k!)^{2}} T_{k}(x)
\end{equation}
where
\begin{equation}\label{lmda}\lambda_{k}=M q_{k}+N r_{k}+M N s_{k}.\end{equation}
\subsection{Characterization using Bernstein basis}
The Bernstein polynomials have been studied thoroughly and there exist many great enduring works on theses polynomials \cite{Farouki5}. 
They are known for their analytic and geometric properties \cite{Farin,Hoschek}, where the basis are known to be optimally stable. They are all non-negative, $B_{i}^{n}(x)\geq 0,$ $x\in[0,1],$ form a partition of unity (normalization) $\sum_{k=0}^{n}B_{k}^{n}(x)=1,$ satisfy symmetry relation $B_{i}^{n}(x)=B_{n-i}^{n}(1-x),$ have a single unique maximum of 
$\binom{n}{i}i^{i}n^{-n}(n-i)^{n-i}$ at $x=\frac{i}{n},$ $i=0,\dots,n,$ and
their roots are $x=0,1$ with multiplicities. The Bernstein polynomials of degree $n$ can be defined by combining two Bernstein polynomials
of degree $n-1.$ That is, the $k$th $n$th-degree Bernstein polynomial defined by the following recurrence relation
 $B_{k}^{n}(x)=(1-x)B_{k}^{n-1}(x)+xB_{k-1}^{n-1}(x),\hspace{.1in} k=0,\dots,n; n\geq 1$ where $B_{0}^{0}(x)=0$ and $B_{k}^{n}(x)=0$ for $k<0$ or $k>n.$ For more details, see Farouki \cite{Farouki5}.

In addition, it is possible to write Bernstein polynomial of degree $r$ where $r\leq n$ in terms of Bernstein polynomials of degree $n$ using
 the following degree elevation \cite{Farouki3}:
\begin{equation}\label{ber-elv} B_{k}^{r}(x)=\sum_{i=k}^{n-r+k}\frac{\binom{r}{k}\binom{n-r}{i-k}}{\binom{n}{i}}B_{i}^{n}(x),\hspace{.1in}k=0,1,\dots,r.
\end{equation}

Now, to write a generalized Chebyshev-type polynomial $\mathscr{T}_{r}^{(M,N)}(x)$ of degree $r$ as a linear combination of
the Bernstein polynomial basis $B_{i}^{r}(x), i=0,1,\dots,r$ of degree $r$ in explicit closed form, we begin with substituting 
\eqref{binomial} into \eqref{gen-Chebyshev-I-ccomb} to get
\begin{equation*}
\begin{aligned}
\mathscr{T}_{r}^{(M,N)}(x)
&=\frac{(2r-1)!!}{(2r)!!}\sum_{i=0}^{r}(-1)^{r-i}\eta_{i,r}B_{i}^{r}(x)+\sum_{k=0}^{r}
\frac{(2k)!\lambda_{k}}{2^{2k}(k!)^{2}}\sum_{j=0}^{k}(-1)^{k-j}\eta_{j,k}B_{j}^{k}(x).
\end{aligned}
\end{equation*}
where 
\begin{equation*}\eta_{i,r}=\frac{\binom{r-\frac{1}{2}}{i}\binom{r-\frac{1}{2}}{r-i}}{\binom{r}{i}},\hspace{.1in} i=0,1,\dots,r.\end{equation*}
This shows that the generalized Chebyshev-type polynomial $\mathscr{T}_{r}^{(M,N)}(x)$ of degree $r$ can be written in the Bernstein basis form.

Now, by expanding the right-hand side and using \eqref{n-half} with some simplifications, we have
\begin{align*}
\eta_{i,r}
=\frac{\left(r-\frac{1}{2}\right)!\left(r-\frac{1}{2}\right)!}
       {(i-\frac{1}{2})!(r-i-\frac{1}{2})!r!}=\frac{(2r-1)!!(2r-1)!!}
       {2^{r}r!(2i-1)!!(2(r-i)-1)!!}.
\end{align*}
Using the fact $(2n)!=(2n-1)!!2^{n}n!$ we get $\eta_{i,r}=\frac{\binom{2r}{r}\binom{2r}{2i}}{2^{2r}\binom{r}{i}}.$

It is clear that $\eta_{0,r}=\frac{1}{2^{2r}}\binom{2r}{r}.$ With simple combinatorial identities simplifications we have
$$\eta_{i-1,r}=\frac{(i-\frac{1}{2})}{(r-i+\frac{1}{2})}\eta_{i,r}.$$
Thus we have the following theorem.

\begin{theorem}\label{gen-jacinBer form}
For $M,N\geq 0,$ the generalized Chebyshev-type polynomials $\mathscr{T}_{n}^{(M,N)}(x)$ of degree $n$ have the following
Bernstein representation:
\begin{equation*}\label{gen-jac-inBer-r}
\mathscr{T}_{n}^{(M,N)}(x)=\frac{(2n-1)!!}{(2n)!!}\sum_{i=0}^{n}(-1)^{n-i}\eta_{i,n}B_{i}^{n}(x)+\sum_{k=0}^{n}
\frac{(2k)!\lambda_{k}}{2^{2k}(k!)^{2}}\sum_{j=0}^{k}(-1)^{k-j}\eta_{j,k}B_{j}^{k}(x)\end{equation*}
where
$\lambda_{k}=M q_{k}+N r_{k}+M N  s_{k}$
and
$\eta_{i,n}=\frac{\binom{2n}{n}\binom{2n}{2i}}{2^{2n}\binom{n}{i}}, i=0,1,\dots,n$ where
$\eta_{0,n}=\frac{1}{2^{2n}}\binom{2n}{n}.$
Moreover, the coefficients $\eta_{i,n}$ satisfy the recurrence relation
\begin{equation}\label{recur}\eta_{i,n}=\frac{(2n-2i+1)}{(2i-1)}\eta_{i-1,n}, \hspace{.1in}i=1,\dots,n.\end{equation}
\end{theorem}

It is worth mentioning that Bernstein polynomials can be differentiated and integrated easily as 
$$\frac{d}{dx}B_{k}^{n}(x)=n[B_{k-1}^{n-1}(x)-B_{k}^{n-1}(x)],\hspace{.02in}n\geq 1,\hspace{.05in}\text{and}\hspace{.05in}
\int_{0}^{1}B_{k}^{n}(x)dx=\frac{1}{n+1},\hspace{.02in} k=0,1,\dots,n.$$

Rababah \cite{rababah4} provided some results concerning integrals of univariate Chebyshev-I and Bernstein polynomials. In the following we consider
 integration of the weighted generalized Chebyshev-type with Bernstein polynomials, 
$$I=\int_{0}^{1}x^{-\frac{1}{2}}(1-x)^{-\frac{1}{2}}B_{r}^{n}(x)\mathscr{T}_{i}^{(M,N)}(x)dx.$$
By using \eqref{gen-jac-inBer-r}, the integral can be simplified to
\begin{equation*}\begin{aligned}
I&=\int_{0}^{1}(1-x)^{-\frac{1}{2}}x^{-\frac{1}{2}}\binom{n}{r}x^{r}(1-x)^{n-r}
\frac{(2i)!}{2^{2i}(i!)^{2}}
\sum_{k=0}^{i}(-1)^{i-k}\frac{\binom{i-\frac{1}{2}}{k}\binom{i-\frac{1}{2}}{i-k}}{\binom{i}{k}}B_{k}^{i}(x)\\
&+\sum_{d=0}^{i}\lambda_{d}\int_{0}^{1}(1-x)^{n-r-\frac{1}{2}}x^{r-\frac{1}{2}}\binom{n}{r}
\frac{(2d)!}{2^{2d}(d!)^{2}}
\sum_{j=0}^{d}(-1)^{d-j}\frac{\binom{d-\frac{1}{2}}{j}\binom{d-\frac{1}{2}}{d-j}}{\binom{d}{j}}B_{j}^{d}(x)dx,
\end{aligned}
\end{equation*}
where $\lambda_{d}$ defined in \eqref{lmda}. By reordering the terms we get
\begin{equation*}\begin{aligned}
I&=\binom{n}{r}\frac{(2i)!}{2^{2i}(i!)^{2}}
\sum_{k=0}^{i}(-1)^{i-k}\binom{i-\frac{1}{2}}{k}\binom{i-\frac{1}{2}}{i-k}\int_{0}^{1}x^{r+k-\frac{1}{2}}(1-x)^{n+i-r-k-\frac{1}{2}}dx\\
&+\sum_{d=0}^{i}\binom{n}{r}\frac{(2d)!\lambda_{d}}{2^{2d}(d!)^{2}}
\sum_{j=0}^{d}(-1)^{d-j}\binom{d-\frac{1}{2}}{j}\binom{d-\frac{1}{2}}{d-j}\int_{0}^{1}x^{r+j-\frac{1}{2}}(1-x)^{n+d-r-j-\frac{1}{2}}dx.
\end{aligned}
\end{equation*}
The integrals in the last equation are the Beta functions $B(x_{i},y_{i})$ with $x_{1}=r+k+\frac{1}{2},$ $y_{1}=n+i-r-k+\frac{1}{2},$
$x_{2}=r+j+\frac{1}{2},$ and $y_{2}=n+d-r-j+\frac{1}{2}.$

Hence, the following theorem provides a closed form of the integration of the weighted generalized Chebyshev-type with respect to the Bernstein polynomials.
\begin{theorem}\label{gen-int-ber-jac}
Let $B_{r}^{n}(x)$ be the Bernstein polynomial of degree $n$ and $\mathscr{T}_{i}^{(M,N)}(x)$ be the generalized Chebyshev-type polynomial of degree $i,$ then
for $i,r=0,1,\dots,n$ we have
\begin{equation*}
\begin{aligned}
&\int_{0}^{1}(x-x^{2})^{-\frac{1}{2}}B_{r}^{n}(x)\mathscr{T}_{i}^{(M,N)}(x)dx\\
&=\binom{n}{r}\frac{(2i)!}{2^{2i}(i!)^{2}}\sum_{k=0}^{i}(-1)^{i-k}\binom{i-\frac{1}{2}}{k}\binom{i-\frac{1}{2}}{i-k}
B(r+k+\frac{1}{2},n+i-r-k+\frac{1}{2})\\
&+\sum_{d=0}^{i}\binom{n}{r}\frac{(2d)!\lambda_{d}}{2^{2d}(d!)^{2}}\sum_{j=0}^{d}(-1)^{d-j}
\binom{d-\frac{1}{2}}{j}\binom{d-\frac{1}{2}}{d-j}
B(r+j+\frac{1}{2},n+d-r-j+\frac{1}{2})
\end{aligned}
\end{equation*}
where $\lambda_{k}=M q_{k}+N r_{k}+M N s_{k}$ and $B(x,y)$ is the Beta function.
\end{theorem}

\section{Applications}
The analytic and geometric properties of the Bernstein polynomials made them important for the development of B\'{e}zier curves and
surfaces in Computer Aided Geometric Design. The Bernstein polynomials are the standard basis for the B\'{e}zier representations of curves
and surfaces in CAGD. However, the Bernstein polynomials are not orthogonal and could not be used effectively in the least-squares approximation \cite{Rice},
and thus the calculations performed in obtaining the least-square approximation polynomial of degree $n$ do not reduce the calculations to
 obtain the least-squares approximation polynomial of degree $n+1.$ So, the method of least squares approximation accompanied by orthogonal polynomials has 
been introduced. 
\begin{definition}
For a function $f(x),$ continuous on $[0,1]$ the least square approximation requires finding a polynomial (Least-Squares Polynomial)
$$p_{n}(x)=a_{0}\varphi_{0}(x)+a_{1}\varphi_{1}(x)+\dots+a_{n}\varphi_{n}(x)$$ 
that minimize the error
$$E(a_{0},a_{1},\dots,a_{n})=\int_{0}^{1}[f(x)-p_{n}(x)]^{2}dx.$$
\end{definition}

For minimization, the partial derivatives must satisfy $\frac{\partial E}{\partial a_{i}}=0, i=0,\dots,n.$
These conditions give rise to a system of $(n+1)$ normal equations in $(n+1)$ unknowns: $a_{0},a_{1},\dots,a_{n}.$
Solution of these equations will yield the unknowns of the least-squares polynomial $p_{n}(x).$ 
By choosing $\varphi_{i}(x)=x^{i},$  then the coefficients of the normal equations give the Hilbert matrix
that has round-off error difficulties and notoriously ill-condition for even modest values of $n.$
 However, choosing $\{\varphi_{0}(x),\varphi_{1}(x),\dots,\varphi_{n}(x)\}$ to be orthogonal simplifies the least-squares approximation problem,  the matrix of the normal equations will be diagonal, 
which make the numerical calculations more efficient. See \cite{Rice} for more details on the least squares approximations.

{\bf Acknowledgements.}
The author would like to thank the referees for their valuable comments which helped to improve the manuscript.

\end{document}